\documentclass[11pt]{amsart}
\usepackage{amsthm,enumerate}
\usepackage{amssymb}
\usepackage{color}
\usepackage{mathtools}
\usepackage[bookmarks, bookmarksdepth=2, colorlinks=true, linkcolor=blue, citecolor=blue, urlcolor=blue]{hyperref}

\usepackage[utf8]{inputenc}
\usepackage[T1]{fontenc}
\usepackage{lmodern}
\usepackage{cleveref}
\usepackage{hyperref}
\usepackage{algorithm}
\usepackage{comment}
\usepackage{url, mathrsfs}
\usepackage{amsmath}
\usepackage{graphicx, 
epstopdf}
\AppendGraphicsExtensions{.gif}
\usepackage{subfig}
\usepackage{color}
\usepackage{bbm}
\usepackage{subfloat}
\usepackage[draft]{fixme}
\usepackage{bm} 
\usepackage{bbm} 
\usepackage{epigraph}
\usepackage{endnotes}
\usepackage{ifpdf}
\usepackage{array}
\usepackage{wrapfig}
\usepackage{longtable}
\usepackage{multicol}

\usepackage{multicol}
\usepackage{multirow}
\usepackage{graphicx}
\usepackage{epstopdf}
\usepackage{tikz}
\usepackage{tikz-cd}
\usepackage{makecell}
\usepackage{fullpage}
\usepackage{media9}%
\theoremstyle{plain}
\newtheorem{theorem}{Theorem}[section]
\newtheorem{proposition}[theorem]{Proposition}
\newtheorem{corollary}[theorem]{Corollary}
\newtheorem{lemma}[theorem]{Lemma}

\newtheorem*{claim*}{Claim}

\newtheoremstyle{case}{}{}{}{}{}{:}{ }{}
\theoremstyle{case}

\theoremstyle{definition}
\newtheorem{definition}[theorem]{Definition}
\newtheorem{example}[theorem]{Example}
\newtheorem{remark}[theorem]{Remark}

\newcommand \R{{\mathbb R}}
\newcommand \C{{\mathbb C}}

\newcommand \K{{\mathcal K}}

\def \Re{\mathsf{Re}}

\def \even{\mathrm{even}}
\def \odd{\mathrm{odd}}

\DeclareMathOperator{\SL}{\mathrm{SL}}

\DeclareMathOperator{\inter}{int}

\newcommand \0{{\mathbf{0}}}

\def \H{{\mathcal{H}}}
\def \Hur{{\mathrm{Hur}}}

\title{$\K$-Lorentzian and $\K$-CLC Polynomials in Stability Analysis}

\author{Papri Dey}
\address{Applied Mathematics, Baskin School of Engineering, University of California, Santa Cruz}
\date{}

\begin{document}

\begin{abstract}
We investigate the class of $\K$-Lorentzian polynomials, which generalizes the well-known class of Lorentzian polynomials. As demonstrated in \cite{GPlorentzian}, these polynomials are equivalent to $\K$-completely log-concave ($\K$-CLC) forms. Exploiting the structure of directional derivatives and Hessians, we prove that the Hessian
matrices of $\K$-CLC forms satisfy a generalized Perron--Frobenius property over a self-dual cone $\K$: they admit a distinguished eigenvector in $\inter\K$. We introduce a line-restriction perspective on $\K$-CLC. This leads to a hierarchy of log-concavity and Newton-type log-concavity inequalities for the coefficients of univariate restrictions. This provides useful necessary conditions for $\K$-CLC in terms of directional derivatives.

As a key application, we derive degree-dependent criteria for Hurwitz stability of $\K$-CLC univariate restrictions. In particular, we prove that all strictly $\K$-Lorentzian polynomials of degree $d \le 4$ are Hurwitz-stable over $\K$, and we obtain refined inequalities that give sufficient conditions for Hurwitz stability in degrees $d \ge 5$.
These results recover and extend known stability criteria in the classical Lorentzian case $\K = \R^n_+$ to the general conic setting. Moreover, these Hurwitz–stability and Newton-type log-concavity coefficient criteria yield simple sufficient conditions for asymptotic stability of linear time-invariant systems via their characteristic polynomials.

\end{abstract}

\maketitle

\section{Introduction}

$\K$-Lorentzian polynomials were introduced and studied by Brändén and Leake in \cite{Leake}, and by Blekherman and Dey in \cite{GPlorentzian}. The approach in \cite{GPlorentzian} not only generalizes Lorentzian polynomials introduced by Huh and Brändén \cite{Brandenlorentzian} and, independently, by Anari, Oveis Gharan and Vinzant under the name of completely log-concave polynomials \cite{Cynthialog}, but also provides a unified framework for extending this class to arbitrary proper convex cones.

In this work we use the term \emph{$\K$-Lorentzian} for the homogeneous setting and \emph{$\K$-completely log-concave} ($\K$-CLC) for the possibly non-homogeneous setting. A central result of \cite{GPlorentzian} is that, for any proper cone $\K$, the class of $\K$-Lorentzian forms coincides with the class of $\K$-CLC forms of the same degree; see \Cref{them:plsclc}. This gives a powerful bridge between convex-analytic log-concavity on cones and the algebraic signature conditions that define Lorentzian behavior via iterated directional derivatives and quadratic forms.

Our first goal in this part is to refine the structural picture of $\K$-Lorentzian and $\K$-CLC polynomials. For homogeneous forms, we show that the Hessian matrices evaluated on $\K$ have a very rigid eigenvalue structure: they have exactly one positive eigenvalue, and, over self-dual cones, satisfy a generalized Perron–Frobenius property (\Cref{them:nonsingHes}, \Cref{them:clcbdd}, \Cref{cor:gpf}). This connects $\K$-Lorentzian polynomials to $\K$-nonnegative and $\K$-irreducible matrices, and to the cone-theoretic Perron–Frobenius theory developed in \cite{Kreinlinear,Barkeralgebraic}.

Next, we introduce a linewise notion of $\K$-complete log-concavity: for each $x \in \inter \K$ and $v \in \K$, we consider the univariate restriction
\[
  t \longmapsto f(x+tv),
\]
and require all iterated derivatives in $t$ to be log-concave on the appropriate interval. We show that $\K$-CLC in the sense of \Cref{def:clc} implies linewise $\K$-CLC (\Cref{prop:kclcimplieslinewisekclc}). Combined with sharp univariate characterizations of complete log-concavity, this yields log-concavity of the factorial-weighted coefficient sequences of $f(x+tv)$ and, consequently, log-concavity of the directional derivative sequence
\[
  \bigl\{ f(x), D_v f(x), D_v^2 f(x), \dots, D_v^d f(x) \bigr\}
\]
for any $x,v \in \inter \K$ (\Cref{thm:logconcave}). We also obtain positivity of the coefficients of these univariate restrictions for strictly $\K$-CLC polynomials (\Cref{prop:noncoeffs}, \Cref{cor:noncoeffs}).

Then we connect these structural inequalities to classical Hurwitz stability of univariate polynomials. Using Routh–Hurwitz and Li\'enard–Chipart criteria together with the log-concavity relations among coefficients, we show that strict $\K$-CLC implies Hurwitz stability in low degrees. Specifically, we prove that every univariate strictly CLC polynomial of degree $d \leq 4$ with positive coefficients is Hurwitz-stable (\Cref{prop:Hurwitzlowerdegree}), and we obtain refined sufficient conditions for Hurwitz stability in degree $d = 5$ (\Cref{them:Hurwitz}) based on strengthened log-concavity inequalities.

In the final part of this paper our Newton-type
coefficient inequalities for strictly CLC and strictly $\K$-Lorentzian
polynomials yield simple sufficient tests ensuring stability of a linear time-invariant system (LTI) directly from the coefficients of $\chi_A(t)$, (\Cref{them:lconstable} and \Cref{them:lconstable5}). 

\smallskip

\noindent\textbf{Organization.} The rest of this paper is structured as follows. In \Cref{sec:sota} we recall basic definitions on proper cones, $\K$-nonnegativity, $\K$-irreducibility, and the equivalence between $\K$-Lorentzian and $\K$-CLC forms. In \Cref{sec:prop} we develop several key properties of $\K$-Lorentzian polynomials, including Hessian characterizations, cone-theoretic Perron–Frobenius conditions, and a line-restriction framework for $\K$-CLC, leading to log-concavity of directional derivatives. In \Cref{sec:hpp} we relate these inequalities to the classical half-plane property and Hurwitz stability for univariate restrictions, deriving degree-dependent criteria for stability. \Cref{sec:dynamics} applies the Hurwitz–stability criteria and Newton-type log-concavity inequalities developed earlier to the characteristic polynomials of linear time-invariant (LTI) systems, yielding simple sufficient conditions for (asymptotic) stability expressed directly in terms of their coefficients.

\section{State of the Art} \label{sec:sota}
A nonempty convex set $\K \subseteq \R^{n}$ is said to be a cone if $c \K \subseteq \K$ for all $c \geq 0$. A cone $\K$ is called proper if it is  closed (in the Euclidean topology on $\R^{n}$), pointed (i.e.,
$\K \cap (-\K) = \{0\}$), and solid (i.e., the topological interior of $\K$, denoted as $\inter \K$, is nonempty). 

Let $\R[x]$ represent the space of $n$-variate polynomials over $\R$, and $\R[x]_{\leq d}$ represent the space of $n$-variate polynomials over $\R$ with degree at most $d$ and $\R[x]^d_n$ denote the set of real homogeneous polynomials (aka forms) in $n$ variables of degree $d$. 
A polynomial $f \in \R[x]$ is said to be (strictly) \emph{log-concave} at a point $a\in \R^n$ if $f(a)>0$, and $\log f$ is a (strictly) concave function at $a$, i.e. the Hessian of $\log f$ is negative semidefinite (negative definite) at $a$. A polynomial $f \in \R[x]$ is log-concave on a proper cone $\K\subset \R^n$ if $f$ is log-concave at every point of the interior of $\K$. By convention, the zero polynomial is log-concave (but not strictly log-concave) at all points of $\R^n$.

If $\K = \K^{\ast}$, a self-dual cone, then $K$ is closed, pointed, and full dimensional\cite{Barkerselfdual}. We denote quadratic forms by a lower case letter, and the matrix of the quadratic form by the corresponding upper case letter, i.e. if $q$ is a quadratic form, then its matrix is $Q$ and $q(x)=x^tQx$.

For a point $a\in \R^n$ and $f\in \R[x]$, $D_a f$ denotes the directional derivative of $f$ in direction $a$: $D_a f=\sum_{i=1}^n a_i\frac{\partial f}{\partial x_i}$. Here are the definitions of the $\K$-completely log-concave polynomials and $\K$-Lorentzian forms.
\begin{definition} \cite{GPlorentzian} \label{def:clc}
A polynomial (form) $f \in \R[x]_{\leq d}$ is called a $\K$-completely log-concave aka $\K$-CLC (form) on a proper convex cone $\K$ if for any choice of $a_{1}, \dots, a_{m} \in \K$, with $m\leq d$, we have that $D_{a_{1}} \dots D_{a_m} f$ is log-concave on $\inter \K$.  
A polynomial (form) $f \in \R[x]$ is strictly $\K$-CLC if for any choice of $a_{1}, \dots, a_{m} \in \K$, with $m\leq d$,  $D_{a_{1}} \dots D_{a_m} f$ is strictly log-concave on all points of $\K$. 
\end{definition}
\begin{definition} \cite{GPlorentzian} \label{def:pls}
Let $\K$ be a proper convex cone. A form $f \in \R[x]_n^d$ of degree $d\geq 2$ is said to be $\K$-Lorentzian if for any $a_{1},\dots,a_{d-2} \in  \inter \K$,  the quadratic form $q=D_{a_{1}} \dots D_{a_{d-2}}f$ satisfies the following conditions:  
\begin{enumerate}
 \item  The matrix $Q$ of $q$ has exactly one positive eigenvalue.
\item  For any $x,y \in\inter \K$ we have $y^t Q x=\langle y,Qx\rangle>0$.
\end{enumerate}
For degree $d \leq 1$ a form is $\K$-Lorentzian if it is nonnegative on $\K$.  
\end{definition} 
\begin{definition} \label{def:k-irre}
Let $\K$ be a proper convex cone in $\R^{n}$. 
An $n\times n$ matrix $A$ is $\K$-nonnegative if $Ax \in \K$ for any $x \in \K$, i.e. $A(\K)\subseteq \K$. A matrix $A$ is $\K$-positive if $Ax \in \inter \K$ for any nonzero $x \in \K$. 
A $\K$-nonnegative matrix $A$ is called $\K$-irreducible if $A$ leaves no (proper) face
of $\K$ invariant. A $\K$-nonnegative matrix which is not $\K$-irreducible is called $\K$-reducible, cf. \cite{gantmacher}, \cite{Vandergraftspectral}.
\end{definition}
\begin{lemma}\label{lem:quadratic} \cite[Lemma 3.1]{GPlorentzian}
A quadratic form $q(x)=x^{t}Qx\in \R[x]_n^2$ is $\K$-Lorentzian if and only if the matrix $Q$ of $q$ has exactly one positive eigenvalue and one of the following conditions holds:
\begin{enumerate}[(a)]
\item for all $x\in \K$ we have $Qx\in \K^*$, i.e. $Q(\K)\subseteq{\K^*}$.
     \item $y^tQx \geq 0$ for all $x,y$ spanning extreme rays of $\K$.
    \item $x^tQx > 0$ for all $x \in \inter \K$.
\end{enumerate}
\end{lemma}
Similar to copositive matrices over the nonnegative orthant, we define $\K$-copositive matrices over $\K$, see \cite{stabilitylevi}. A matrix $Q$ is (strictly) $\K$-copositive if $x^{t}Qx (>0) \geq 0$ for all $x \in \K$.
\begin{corollary} \label{cor:copositive}
If a quadratic form $q(x)=x^{t}Qx\in \R[x]_n^2$ is (strictly) $\K$-Lorentzian, then $Q$ is (strictly) $\K$-copositive matrix.
\end{corollary}

For a self-dual cone $\K$ we have the following result.
\begin{corollary}\label{cor:kpos} \cite[Cor 3.8]{GPlorentzian}
Let $\K$ be a self-dual cone. Then the quadratic form $q=x^{t}Qx$ is $\K$-Lorentzian if and only if the matrix $Q$ of $q$ has exactly one positive eigenvalue and $Q$ is $\K$-nonnegative.
\end{corollary}

\begin{proposition} \label{cor:qclcirre} \cite[Cor 3:10]{GPlorentzian}
Let $\K\subset \R^n$ be a self-dual cone. Then $q=x^{t}Qx$ is $\K$-Lorentzian if and only if the matrix $Q$ of $q$ has exactly one positive eigenvalue and $Q$ is either nonsingular and $\K$-irreducible, or singular and $\K$-nonnegative.
\end{proposition} 
\begin{remark} \label{remark:quadlorentzian}
Then $q$ is a quadratic Lorentzian (aka stable) polynomial if and only if the matrix $Q$ of $q$ has exactly one positive eigenvalue and $Q$ is nonnegative matrix. 
\end{remark}
Furthermore, we use the following equivalence result between $\K$-Lorentzian and $\K$-CLC polynomials, proved in \cite[Theorem 4:10]{GPlorentzian}.
\begin{theorem} \label{them:plsclc}
Let $\K\subset \R^n$ be a proper cone. A form $f\in \R[x]_n^d$ is $\K$-Lorentzian if and only if $f$ is $\K$-CLC. Equivalently, $\SL^d(\K)=\mathrm{CLC}^d (\K)$.
\end{theorem}
\section{Properties of \texorpdfstring{$\K$}{}-Lorentzain polynomials} \label{sec:prop}
\subsection{Proper Convex Cones}
Using \Cref{lem:quadratic} we show that $\K$-nonnegativity is a sufficient condition for $\K$-Lorentzian polynomials if $\K$ is contained in its dual cone, $\K^{\ast}$.
\begin{corollary} \label{cor:conecontainsuff}
Consider the quadratic form $q(x)=x^{t}Qx \in \R[x]_n^2$ and $\K$ be a proper cone such that $\K \subseteq \K^{\ast}$. If the matrix $Q$ of $q$ has exactly one positive eigenvalue and $Q$ is $\K$-nonnegative, then $q$ is a $\K$-Lorentzian on the proper convex cone $\K$.
\end{corollary}
\begin{proof}
    Since $Q$ is $\K$-nonnegative, so $Q(\K) \subseteq \K \subseteq \K^{\ast}$. Then the rest follows from \Cref{lem:quadratic}.
\end{proof}

Based on \Cref{lem:quadratic} we show that $\K$-nonnegativity is a necessary condition for $\K$-Lorentzian polynomials if $\K$ contains its dual cone, $\K^{\ast}$.
\begin{corollary} \label{cor:conecontain}
If a quadratic form $q(x)\in \R[x]_n^2$ is $\K$-Lorentzian on a proper convex cone $\K$ and $\K^{\ast} \subseteq \K$, the matrix $Q$ of $q$ has exactly one positive eigenvalue and $Q$ is $\K$-nonnegative.
\end{corollary}
\begin{proof} By \Cref{lem:quadratic} we have the matrix $Q$ of $q$ has exactly one positive eigenvalue and $Q(\K)\subseteq{\K^*}$. Thus, $Q(\K)\subseteq \K^* \subseteq \K$. Therefore, $Q$ is $\K$-nonnegative by \Cref{def:k-irre}.  
\end{proof}
\begin{remark}
   The converse to \Cref{cor:conecontain} is not true. It's easy to construct counterexamples where $\K^{\ast} \subseteq \K$ and $Q(\K) \subseteq \K$ but $Q(\K) \not \subseteq \K^{\ast}$. %see notebook 1, June, 2024
\end{remark}
It's shown that any closed convex cone contains its dual cone, i.e., $\K^{\ast} \subset \K$, if and only if there exist $y,z \in \K$ such that $x=y-z$ and $ \langle y, z \rangle =0$, cf. \cite{Remarkcopositive} for detail. This is used to show that every $x \in \R^{n}$ has an (unique) orthogonal decomposition on $\K$ in \cite{Cross} which have many well-known consequences. Additionally, the sum of two $\K$-nonnegative matrices is $\K$ nonnegative. Furthermore,  if $A$ and $B$ are $\K$-nonnegative and $A$ is $\K$-irreducible, then so is $A + B$, cf. \cite[Chap-1]{Bermannonnegative}. However, the sum of two log-concave polynomials may not be log-concave.
The sufficient condition for the sum of two $\K$-Lorentzian polynomials to also be $\K$-Lorentzian can be found in \cite[Lemma 3.3]{Cynthialog3} for the nonnegative orthant, and \cite[Theorem 4.15]{GPlorentzian} for any proper convex cone $\K$. 
\begin{theorem} \label{them:sum} \cite[Theorem 4.15]{GPlorentzian}
Let $f,g \in \R^{d}_{n}[x]$ be two $\K$-Lorentzian polynomials over the self-dual cone $\K$. Then the sum $f+g$ is $\K$-Lorentzian if there exist vectors $b,c \in \K$ such that $D_bf=D_cg \not \equiv 0$.
\end{theorem}
\subsection{Perron-Frobenius Theorem for Cones} Let $\K$ be a proper convex cone. We can leverage the relationship between $\K$-Lorentzian polynomials and the generalized Perron-Frobenius theorem \cite{Kreinlinear}, \cite{Barkeralgebraic} to distinguish potential $\K$-Lorentzian polynomial candidates from given log-concave polynomials. This also provides valuable insight into the structure of the cone K. 
\begin{theorem}(Generalization of Perron-Frobenius Theorem)  \label{them:genpf} 
  \cite{Kreinlinear}, \cite{Barkeralgebraic}: Let $A$ be a $\K$-irreducible matrix with spectral radius $\rho$. Then
\begin{enumerate}
\item $\rho$ is a simple positive eigenvalue of $A$,
\item There exists a (up to a scalar multiple) unique $\K$-positive (right)
eigenvector $u$ of $A$ corresponding to $\rho$,
\item $u$ is the only $\K$-semipositive eigenvector for $A$ (for any eigenvalue),
\item $\K \cap (\rho I-A)\R^{n} =\{0\}$.
\end{enumerate}
\end{theorem}
By applying \cite[Lemma 4.12]{GPlorentzian} and \cite[Prop:4.13]{GPlorentzian} we derive the following characterization of the Hessian matrices evaluated at any point over the cone $\K$ for any $\K$-Lorentzian polynomials where $\K$ is a proper convex cone.

\begin{theorem} \label{them:nonsingHes}
Let $f \in \R[x]_n^d$ be a nonzero form of degree $d \geq 2$, and $\K$ be a proper convex cone in $\R^{n}$. 
\begin{enumerate}[(a)]
\item If $f$ is a strictly $\K$-CLC form, then $H_{f}(a)$ is nonsingular and has exactly one positive eigenvalue for all $a \in \K$.  Also the quadratic form $x^{t} H_{f}(a)x$ is negative definite on $(H_{f}(a)b)^{\perp}$ for every $a, b \in  \K$ such that $H_{f}(a)b \neq 0$. 
\item If $f$ is a $\K$-CLC form, then $H_{f}(a)$ has exactly one positive eigenvalue for all $a \in \inter \K$. Also the quadratic form $x^{t} H_f(a)x$ is negative semidefinite on $(H_f(a)b)^{\perp}$ for every $a \in \inter \K$, and $b \in \K$ such that $H_f(a)b \neq 0$. 
\end{enumerate}
\end{theorem}
In special cases using the \Cref{def:pls} and \cite[Them:5.1]{GPlorentzian}
we have some necessary conditions for $\K$-Lorentzian polynomials over the self-dual cone $\K$.
\begin{theorem} \label{them:clcbdd}
Let $f \in \R[x]_n^d$ be a $\K$-Lorentzian polynomial of degree $d \geq 2$ over the self-dual cone $\K$. Then for all $a \in  \K$

    \begin{enumerate}[(a)]
        \item $H_{f}(a)$ is either nonsingular and $\K$-irreducible, or singular and $\K$-nonnegative.
        \item $H_{f}(a)$ has exactly one positive eigenvalue.
        \end{enumerate}
        \end{theorem}

By \Cref{them:genpf} we have more stronger necessary conditions for $\K$-Lorentzian polynomials over the self-dual cone $\K$.
\begin{corollary} \label{cor:gpf}
Let $f \in \R^{d}_{n}[x]$ be a $\K$-Lorentzian polynomial over the self-dual cone $\K$. Then for any $a \in \K$, the nonsingular Hessian $H_f(a)$ satisfies generalized Perron-Frobenius property. 
\end{corollary}
\begin{proof} Since $f$ is a $\K$-Lorentzian polynomial, so for any $a \in  \K$, the nonsingular Hessian $H_{f}(a)$ is $\K$-irreducible matrix and has exactly one positive eigenvalue.
Therefore, if $H_{f}(a)$ is nonsingular, and the positive eigenvalue corresponds to the spectral radius, then by \Cref{them:genpf}, its corresponding eigenvector must lie in the interior of the cone $\K$.
\end{proof}

The geometry of the set of symmetric Perron–Frobenius matrices has been studied in \cite{SymmetricPFstructure}. In fact, it is not convex but is star convex with the identity matrix as the center. It's a natural question to describe the geometry of the set of symmetric Perron–Frobenius matrices which have exactly one positive eigenvalue. We leave this question for future work.

\subsection{Properties of \texorpdfstring{$\K$}{}-CLC polynomials}

Here we propose another definition for $\K$-CLC multivariate polynomials by restricting them on certain lines. This allows us to show the log-concavity property among the coefficients of the corresponding univariate restrictions of a $\K$-CLC polynomial. Recall that a univariate polynomial $f(t)$ is log-concave over an interval $I$ if $f(t) > 0$ for all $t \in I$ and $f(t)f^{''}(t)\leq (f^{'}(t))^{2}$ for all $t \in I$. Equivalently, a univariate polynomial $f(t)$ is log-concave over an interval $I$ if $f(t) > 0$ for all $t \in I$ and $\frac{d^{2}}{dt^2} \log f(t) \leq 0$ for all $t \in I$. So, we can define $\K$-CLC for any real multivariate polynomial as in \Cref{def:clc}, need not be homogeneous. If they are homogeneous, the set of $\K$-Lorentzian polynomials of degree $d$ are equivalent to the set of $\K$-CLC forms , cf. \cite{GPlorentzian}).
On the other hand, over the nonnegative orthant, a characterization for bivariate $\K$-CLC form has been established in the literature, providing a resolution to Mason's conjecture.
\begin{theorem} \label{them:bivclc} \cite{Gurvitsmultivariate}, \cite{Cynthialog3}
$f = \sum_{k=0}^{n} c_kx^{n-k}y^{k} \in \R[x,y]$ is completely log- concave if and only if the sequence of nonnegative coefficients, $\{c_0, \dots , c_n\}$ is ultra log-concave, i.e., for every $1 < k < n$, $\left(\frac{c_{k}}{\binom{n}{k}}\right)^2 \geq \frac{c_{k-1}}{\binom{n}{k-1}} \frac{c_{k+1}}{\binom{n}{k+1}}$.
\end{theorem}

\begin{lemma} \label{lemma:univCLC}
Let $f(t)=\sum_{i=0}^{d}a_{i}t^{i}$ be a univariate polynomial. The univariate polynomial $f$ is CLC over $t \geq 0$ if and only if $a_i >0$ and the sequence $\{0!a_0, \dots, d!a_d\}$ is a log-concave sequence, i.e., $ia_{i}^{2} \geq (i+1)a_{i-1}a_{i+1}$ for all $i=0(1)d$. 
\end{lemma}
\proof It's straightforward to verify that $f(t)=\sum_{i=0}^{d}a_{i}t^{i}$ with $a_i >0$ is CLC for $t \geq 0$ if and only if its homogenization $f(t,u)=\frac{u^{d}}{d!}a_0+\frac{u^{d-1}}{(d-1)!}a_1t+\dots+\frac{u^{0}}{0!}a_dt^d$ is CLC over $\R^{2}_{\geq 0}$. Then the rest follows from \Cref{them:bivclc} and 
the identity $a_{i}^{2}(i!)^2 \geq a_{i-1}a_{i+1}(i-1)!(i+1)! \Leftrightarrow ia_{i}^{2} \geq (i+1)a_{i-1}a_{i+1}$.  \qed

\begin{corollary} \label{cor:reverse}
The polynomial $f(t)=\sum_{i=0}^{d}a_it^{d-i}$ with $a_{i} > 0$ is CLC for all $t \geq 0$ if and only $\{d!a_0, \dots,i!a_{d-i}, \dots, 0!a_d\}$ is a log-concave sequence 
\end{corollary}
\begin{proof}
    Due to \Cref{lemma:univCLC} and symmetry, $((d-k)!a_k)^2 \geq (d-k-1)!(d-k+1)! a_{k-1}a_{k+1} \Leftrightarrow (d-k)a_k^2 \geq (d-k+1)a_{k-1}a_{k+1} \Leftrightarrow ia_{d-i}^{2} \geq (i+1)a_{d-i-1}a_{d-i+1}$.
\end{proof}

It is well known that convexity can be characterized in terms of
one-dimensional convexity.  More precisely,  a function $f:\R^{n} \rightarrow \R$ is convex over a convex set $\K$ if and only if for any $x \in \K$ and $v \in \R^{n}$, the function $g:\R \rightarrow \R, g(t)=f(x+tv)$ is convex as function in $t$ for all $t$ such that $x+tv \in \K$ \cite{Boydconvex}. In the same spirit, we introduce a
line-restriction notion for log-concavity and complete log-concavity
over cones.

\begin{definition}[Linewise $\K$-completely log-concave]  \label{def:clcline}
Let $\K\subseteq\R^n$ be a proper convex cone and let
$f\in\R[x_1,\dots,x_n]$.  For $x\in\operatorname{int}\K$ and
$v\in\K\setminus\{0\}$, set
\[
  I_{x,v}\coloneqq\{\,t\in\R : x+tv\in\operatorname{int}\K\,\},
  \qquad
  g_{x,v}(t)\coloneqq f(x+tv),\quad t\in I_{x,v}.
\]
\begin{itemize}
  \item We say that $f$ is \emph{linewise log-concave over $\K$} if,
  for every $x\in\operatorname{int}\K$ and $v\in\K\setminus\{0\}$, the
  univariate function $g_{x,v}$ is log-concave on $I_{x,v}$.
  \item We say that $f$ is \emph{linewise $\K$-completely log-concave}
  (linewise $\K$-CLC) if, for every $x\in\operatorname{int}\K$ and
  $v\in\K\setminus\{0\}$, the univariate polynomial $g_{x,v}$ is
  completely log-concave on $I_{x,v}$ in the univariate sense (i.e.,
  for every $m\ge 0$, the derivative $g_{x,v}^{(m)}$ is log-concave on
  $I_{x,v}$).
\end{itemize}
\end{definition}

We next show that $\K$-complete log-concavity in the sense of
\Cref{def:clc} implies the linewise property
from \Cref{def:clcline}.

\begin{proposition}\label{prop:kclcimplieslinewisekclc}
Let $\K\subseteq\R^n$ be a proper convex cone, and suppose that
$f\in\R[x_1,\dots,x_n]$ is $\K$-completely log-concave in the sense of
\Cref{def:clc}. Then $f$ is linewise $\K$-completely
log-concave.
\end{proposition}

\begin{proof}
Fix $x\in\operatorname{int}\K$ and $v\in\K\setminus\{0\}$, and let
$I_{x,v}$ and $g_{x,v}$ be as in
\Cref{def:clcline}.  We must show that $g_{x,v}$ is
completely log-concave on $I_{x,v}$ in the univariate sense.

By \Cref{def:clc}, for every integer $m\ge 0$ the
directional derivative
\[
  D_v^m f \;=\; \underbrace{D_v\cdots D_v}_{m\ \text{times}} f
\]
is log-concave on $\operatorname{int}\K$, since each direction $v$
belongs to $\K$.  For $t\in I_{x,v}$ we have $x+tv\in\operatorname{int}\K$, and
\[
  \frac{d^m}{dt^m} g_{x,v}(t)
  \;=\; \frac{d^m}{dt^m} f(x+tv)
  \;=\; D_v^m f(x+tv).
\]
Thus the $m$-th derivative of $g_{x,v}$ is the restriction of the
log-concave function $D_v^m f$ (defined on $\operatorname{int}\K$) to
the line segment $\{x+tv : t\in I_{x,v}\}$.  Since restrictions of a
log-concave function to a line are log-concave on their domain, it
follows that $g_{x,v}^{(m)}$ is log-concave on $I_{x,v}$ for every
$m\ge 0$.

This is precisely the definition of complete log-concavity for the
univariate function $g_{x,v}$ on $I_{x,v}$, and therefore $f$ is
linewise $\K$-completely log-concave.
\end{proof}

Next we provide two necessary conditions for a $f$ to be $\K$-CLC along $v \in \inter \K$.
\begin{proposition} \label{prop:noncoeffs}
   Let $f(x)$ be a nonzero $\K$-CLC over a proper convex cone $\K$. Then for any $x,v \in \inter \K$,  the coefficients of $f(x+tv)$ are positive. 
   \end{proposition}
   \begin{proof} Since $f$ is CLC over $\K$, its directional derivatives are log-concave, i.e., they satisfy the positivity condition at each degree level $d \geq 1$. Therefore, in particular, $D_{v}f(x+tv)=f^{'}(x+tv) >0,$ and $D_{v}D_{v}f(x+tv)=f^{''}(x+tv) >0$ for all $v \in \inter \K$. That enforces all the coefficients of $f(x+tv)=f(x)+t D_v f(x) +\frac{t^{2}}{2} D_{v}^{2}f(x) +\dots+\frac{t^{d}}{d!}D_{v}^{d}f(x)$ to be positive. 
 \end{proof}  
 \begin{corollary}\label{cor:noncoeffs}
Let $f(x)$ be a strictly $\K$-CLC over a proper convex cone $\K$. Then for any $x,v \in \K$, the coefficients of $f(x+tv)$ are positive. In particular, if $f(x)$ is a strictly $\K$-Lorentzian polynomial over a proper convex cone $\K$, then for any $x,v \in \K$, the coefficients of $f(x+tv)$ are positive. 
\end{corollary}
\begin{proof}
    It follows from \Cref{def:clc} and \Cref{prop:noncoeffs}.
\end{proof}

\begin{theorem}\label{thm:logconcave}
Let $K \subseteq \R^n$ be a proper convex cone and let $f(x)$ be a nonzero $K$-CLC polynomial.
Let $d := \deg f$ and let $x, v \in \operatorname{int} K$. Then the sequence
\[
  \bigl\{ f(x),\, D_v f(x),\, D_v^{2} f(x),\, \dots,\, D_v^{d} f(x) \bigr\}
\]
of positive real numbers is log-concave.
\end{theorem}

\begin{proof}
Since $f$ is $K$-CLC in the sense of \Cref{def:clc}, it follows from
\Cref{prop:kclcimplieslinewisekclc} that $f$ is linewise $K$-CLC in the sense of \Cref{def:clcline}. Hence, for any $x, v \in \operatorname{int}K$, the univariate restriction
\[
  g(t) \coloneqq f(x + tv)
\]
is completely log-concave on $t \ge 0$.

By Proposition~\ref{prop:noncoeffs}, for any $x, v \in \operatorname{int}K$ all coefficients of
$g(t)$ are positive. Writing
\[
  g(t) = f(x + tv)
       = f(x) + t D_v f(x) + \frac{t^2}{2!} D_v^2 f(x) + \cdots + \frac{t^d}{d!} D_v^d f(x),
\]
we see that the coefficients are
\[
  a_k = \frac{1}{k!} D_v^k f(x), \qquad k = 0,1,\dots,d.
\]

Since $g$ is univariate CLC with positive coefficients, Lemma~\ref{lemma:univCLC} applies and
implies that the sequence
\[
  \{ 0! a_0,\, 1! a_1,\, \dots,\, d! a_d \}
\]
is log-concave. Substituting $a_k = \frac{1}{k!} D_v^k f(x)$, we obtain
\[
  0! a_0 = f(x),\quad 1! a_1 = D_v f(x),\quad \dots,\quad d! a_d = D_v^d f(x),
\]
so the sequence
\[
  \bigl\{ f(x),\, D_v f(x),\, D_v^{2} f(x),\, \dots,\, D_v^{d} f(x) \bigr\}
\]
is log-concave, as claimed.
\end{proof}

\begin{remark}
The converse of \Cref{thm:logconcave} does not hold.  Even if, for every
$x \in \inter K$ and every $v \in \inter K$, the sequence
\[
  \bigl\{ f(x), D_v f(x), \dots, D_v^d f(x) \bigr\}
\]
is log-concave, this does not imply that $f$ is $K$-CLC.  
\end{remark}
\section{Half-Plane Property} \label{sec:hpp}
Let $\H$ be the open left half-plane $\{x \in \C: \Re \ x<0\}$. %and by $\H$ the closed left half-plane $\{x \in \C: \Re \ x 0\}$. 
A real valued $n \times n$ matrix $A$ is stable if all its eigenvalues belong to the open left half plane. A univariate polynomial $g(t)$ is called Hurwitz-stable if all its roots have negative real part.

Note that Hurwitz stability is closed under inversion, so we have
\begin{proposition} \label{prop:clcline}
  $g(t) \in \R[t]$ is Hurwitz-stable of degree $d$ over an interval $[a,b]$ if and only if $t^{d}g(1/t)$ is Hurwitz-stable over $[a,b]$.    
\end{proposition}
\noindent Let $f(t)=\sum_{k=0}^{d}a_{k}t^{k}$. 
Then the corresponding $d \times d$ Hurwitz matrix is given by $\Hur_{f(t)}=\begin{bmatrix} a_1 & a_3 & a_5 & \dots & \dots & \dots &0&0&0 \\ a_0 & a_2 & a_4 &  & & & \vdots & \vdots & \vdots \\ 0& a_1 & a_3 & & & & \vdots & \vdots & \vdots \\ \vdots & a_0 & a_2 & \ddots & & & 0 & \vdots &\vdots \\ \vdots & 0& a_1 & & \ddots & & a_d &\vdots&\vdots \\ \vdots &\vdots &a_0 & & & \ddots &  a_{d-1} &0 &\vdots \\ \vdots & \vdots &0 & & & & a_{d-2}&a_d& \vdots \\ \vdots & \vdots & \vdots & & & & a_{d-3} &a_{d-1}&0\\0&0&0&\dots&\dots&\dots&a_{d-4}&a_{d-2}&a_{d} \end{bmatrix}$, cf., \cite[Chap-7]{stabilityHurwitz}. \\
\begin{theorem} \label{them:Routh-Hurwitz}
Using Routh-Hurwitz criterion, cf., \cite[Chap 4]{Bellmanselected}, the FAE:
\begin{enumerate}[(a)]
\item $f(t)$ is Hurwitz-stable
\item All the leading principal minors of the corresponding Hurwitz matrix $\Hur_{f(t)}$ are positive.
\end{enumerate}
\end{theorem}
\begin{remark} \label{remark:Routh-Hurwitz}
  By expanding the determinant of $\Hur_{f(t)}$ along the elements of its last column gives $\det (\Hur_{f(t)})=a_d \Delta_{d-1}$ where $\Delta_{k}$ represents the $k \times k$ leading principal minor of $\Hur_{f(t)}$.  Since $a_i >0$, to demonstrate $f(t)$ of degree $d$ is Hurwitz-stable, it suffices to verify that $\Delta_{k}$ are positive for all $1 \leq k \leq d-1$. 
\end{remark}
The following result is a refinement of the Routh-Hurwitz criterion above.

\begin{theorem}(Li{\'e}nard-Chipart) Necessary and sufficient conditions for the polynomial $f(t)=\sum_{k=0}^{d}a_kt^k$ with $a_0>0$ to be stable can be given in any one of the following four forms: \cite[pg-221]{gantmacher}
\begin{enumerate}
    \item $a_k>0, a_{k-2}>0,\dots; \Delta_1>0, \Delta_3 >0, \dots$
    \item $a_k>0, a_{k-2}>0,\dots; \Delta_2>0, \Delta_4 >0, \dots$
    \item $a_k>0, a_{k-1}>0, a_{k-3}>0,\dots; \Delta_1>0, \Delta_3 >0, \dots$
    \item $a_k>0, a_{k-1}>0,a_{k-3}>0,\dots; \Delta_2>0, \Delta_4 >0, \dots$
\end{enumerate}
    
\end{theorem}

A polynomial $f \in \R[x_1,\dots,x_n]$ is called \emph{(Hurwitz-)stable} if
$f(z) \neq 0$ for all $z = (z_1,\dots,z_n) \in \C^n$ with
$\Re(z_i) > 0$ for all $1 \le i \le n$. Equivalently, all zeros of $f$
lie in the closed left half-plane in each coordinate. By convention,
we also regard the zero polynomial as stable.

Equivalently, if $f \in \R[x_1,\dots,x_n]$ is homogeneous and nonzero,
then $f$ is Hurwitz-stable if and only if, for all
$x,v \in (\R_{>0})^n$, the univariate polynomial
\[
t \ \longmapsto\ f(x + t v)
\]
is a Hurwitz-stable polynomial in $t$, that is, all of its zeros have
negative real part; see, e.g., \cite[Lemma~6]{Gregor}, \cite[Prop.~5.2]{Oxley}.

\begin{definition} \label{def:hurstable}
Let $\K \subset \R^n$ be a proper convex cone with nonempty interior
$\inter \K$. A polynomial $f \in \R_n^d[x]$ is said to be
\emph{Hurwitz-stable over $\K$} if either $f \equiv 0$, or else for
every $x,v \in \inter \K$ the univariate polynomial
\[
t \ \longmapsto\ f(x + t v)
\]
is Hurwitz-stable in the univariate sense, i.e., all of its zeros
have negative real part.
\end{definition}

\begin{theorem}\cite[Theorem 5.3]{Gulerhyperbolicconvex}
Any hyperbolic polynomial $f \in \R[x]^{d}_{n}$ is Hurwitz-stable over its hyperbolicity cone $\Lambda_{++}(f,e)$ containing $e$ such that $f(e)>0$.
\end{theorem}
\noindent Here we derive Newton-type inequalities for directional derivatives, which we then feed into stability criteria.
\begin{proposition} \label{prop:logconcaveinequalities}
    Let $f(t)=\sum_{k=0}^{n}a_{k}t^{k}$ be a univariate polynomial with positive $a_k$ for $k=0, \dots,n$. If %$f(t)$ is CLC over $t > 0$, then  
    $a_{k}^{2} \geq \frac{k+1}{k}a_{k-1}a_{k+1}$ for all $1 \leq k \leq d-1$,  the following inequalities hold:
  \begin{enumerate}[(a)]
\item For any $1 \leq k \leq d-2$, $a_{k}a_{k+1} \geq \frac{k+2}{k}a_{k-1}a_{k+2} \Leftrightarrow \frac{a_{k}}{a_{k-1}} \geq \frac{k+2}{k} \frac{a_{k+2}}{a_{k+1}} \Leftrightarrow \frac{a_{k+1}}{a_{k-1}} \geq \frac{k+2}{k} \frac{a_{k+2}}{a_{k}}$.
\item For any $1 \leq k \leq d-2$, $ a_{k}^2 \geq \frac{(k+2)(k+1)}{k(k-1)} a_{k-2}a_{k+2}$.
\item For any $1 \leq k \leq d-3$, $a_ka_{k+2} \geq \frac{k+3}{k}a_{k-1}a_{k+3}$, and for any $1 \leq k \leq d-1$, $a_{k}a_{d-1} \geq \frac{d}{k}a_{k-1}a_d$.
\item $a_ia_j \geq \frac{l}{i} a_ka_l$ if $i+j=k+l$ for any $k < i <j <l$. 
  \end{enumerate}
  \end{proposition}
  \begin{proof}
 $(a)$ 
    \begin{equation*} 
    a_k^2 \geq \frac{k+1}{k}a_{k-1}a_{k+1}, a_{k+1}^2 \geq \frac{k+2}{k+1}a_{k}a_{k+2} \Rightarrow a_{k}a_{k+1} \geq  \frac{k+2}{k}a_{k-1}a_{k+2}
    \end{equation*}
     $(b)$ \begin{eqnarray*} 
    &a_k^2 \geq \frac{k+1}{k}a_{k-1}a_{k+1}, a_{k+1}^2 \geq \frac{k+2}{k+1}a_{k}a_{k+2}, a_{k-1}^2 \geq \frac{k}{k-1}a_{k-2}a_{k} \\
    &\Rightarrow a_{k}^2 \geq  \frac{k+1}{k} \sqrt{\frac{k+2}{k+1}}\sqrt{\frac{k}{k-1}}\sqrt{a_{k-2}a_k}\sqrt{a_{k}a_{k+2}}\\
    &\Rightarrow a_{k}^2 \geq \frac{k+1}{k} \frac{k+2}{k-1} a_{k-2}a_{k+2}
    \end{eqnarray*}
  $(c)$   \begin{eqnarray*} 
    &a_k^2 \geq \frac{k+1}{k}a_{k-1}a_{k+1}, a_{k+1}^2 \geq \frac{k+2}{k+1}a_{k}a_{k+2}, a_{k+2}^2 \geq \frac{k+3}{k+2}a_{k+1}a_{k+3} \\
    &\Rightarrow a_{k}^2 a_{k+2}^2\geq  \frac{k+1}{k} \frac{k+3}{k+2}a_{k-1}a_{k+1}^2a_{k+3}\\
    &\Rightarrow a_{k}^2 a_{k+2}^2\geq  \frac{k+1}{k} \frac{k+3}{k+2} \frac{k+2}{k+1}a_{k-1}a_{k}a_{k+2}a_{k+3}\\
    &\Rightarrow a_{k} a_{k+2}\geq \frac{k+3}{k} a_{k-1}a_{k+3}
    \end{eqnarray*}
    By applying the same process iteratively, we arrive at the final inequality of $(c)$. \\
   $(d)$  This result follows from the established pattern above. 
  \end{proof}
\begin{proposition} \label{prop:logineq}
Let $f(t)=\sum_{k=0}^{n}a_{k}t^{k}$ be a univariate polynomial with positive $a_k$ for $k=0, \dots,n$. If %$f(t)$ is CLC over $t > 0$, then  
    $a_{k}^{2} \geq \frac{k+1}{k}a_{k-1}a_{k+1}$ for all $1 \leq k \leq d-1$, then the following inequalities hold:
\begin{enumerate}[(a)] 
   \item Monotonic chains of adjacent ratios: $\frac{a_2}{a_0} > \frac{a_3}{a_1} > \frac{a_4}{a_2} > \frac{a_5}{a_3} >\dots > \frac{a_k}{a_{k-2}} > \frac{a_{k+1}}{a_{k-1}} > \frac{a_{k+2}}{a_{k}} > \dots $. 
   \item Separated-index ratio chains: $\frac{a_1}{a_0} >\frac{a_3}{a_2} > \frac{a_5}{a_4} >\dots \frac{a_{k-1}}{a_{k-2}} > \frac{a_{k+1}}{a_{k}} > \frac{a_{k+3}}{a_{k+2}} > \dots$
   \item Cross-ratio type inequalities: $\frac{a_{k+1}-a_{k}r}{a_{k-1}-a_{k-2} r} < \frac{a_{k+1}}{a_{k-1}} \Leftrightarrow a_ka_{k-1} >a_{k-2}a_{k+1}$, and \ $\frac{a_{k+1}-a_{k+2}r}{a_{k-1}-a_{k}r} > \frac{a_{k+1}}{a_{k-1}} \Leftrightarrow a_ka_{k+1} >a_{k-1}a_{k+2}$ for any positive $r$. 
   \end{enumerate}
\end{proposition}
\begin{proof}
 $(a)$ The chain of inequalities is obtained by combining \Cref{prop:logconcaveinequalities} $(a)$ and $(d)$.   \\
 $(b)$ It follows from $(a)$ by rearranging the position of $a_i$. \\
 $(c)$ These are derived by applying the identity $\frac{a}{b} > \frac{c}{d} \Leftrightarrow ad >bc$.
 \item 
\end{proof}
\noindent Note that strong log-concavity $\Rightarrow$ very rigid ordering for coefficient ratios and related rational expressions. We translate abstract log-concavity into very concrete, ordered ratio inequalities that can be plugged directly into determinant conditions for Hurwitz matrices.
\begin{lemma} \label{lemma:logconcaveineq}
Let $f(x)$ be a nonzero $\K$-CLC over a proper convex cone $\K$ and $v \in \inter \K$. Then the following inequalities are satisfied.
\begin{enumerate}[(a)]
    \item $(D_{v}^{k}f(x))^2 \geq D_{v}^{k-1}f(x)D_{v}^{k+1}f(x)$ for any $1 \leq k \leq d-1$ and $x \in \inter \K$.
    \item For any $0 \leq k \leq d-3$, $D_{v}^{k+1}f(x)D_{v}^{k+2}f(x) \geq D_{v}^{k}f(x) D_{v}^{k+3}f(x) \Leftrightarrow \frac{D_{v}^{k+1}f(x)}{D_{v}^{k}f(x)} \geq \frac{D_{v}^{k+3}f(x)}{D_{v}^{k+2}f(x)}$.
    \item For any $0 \leq k \leq d-4$, $(D_{v}^{k+2}f(x))^{2} \geq D_{v}^{k}f(x) D_{v}^{k+4}f(x)$.
\end{enumerate}
\end{lemma}
\begin{proof}
$(a)$ Since $f$ is $\K$-CLC, so the sequence $\{f(x),D_vf(x), D_{v}^{2}f(x),\dots, D_{v}^{d}f(x)\}$ of positive real numbers forms a log-concave sequence for any $x \in \inter \K$ by \Cref{thm:logconcave}.  \\
$(b)$ and $(c)$ The inequalities hold for any $x \in \inter \K$ by \Cref{prop:logconcaveinequalities}.
\end{proof}
\noindent Recall that two polynomials $f$ and $g$ are Hurwitz-stable if and only if the product $fg$ is Hurwitz-stable. Additionally, the Hadamard (coefficient-wise) product of two Hurwitz stable polynomials is again Hurwitz stable. Here we report the product rule for $\K$-CLC polynomials. 
 \begin{proposition} \cite[Corollary 8:13]{Hereditary}
 If $f \in \R[x]$ and $g \in \R[x]$ are $\K$-CLC multivariate polynomials, the product $fg$ is $\K$-CLC polynomial
 \end{proposition}
 \noindent The converse of the above statement need not be true in general.
 \begin{example}
Although $f(t)=3+7t+7t^2+4t^3=(1+t+t^2)(3+4t)$ is CLC on $t \geq 0$ but its factor $1+t+t^2$ is not CLC over $t \geq 0$.
 \end{example}
We establish a connection between $\K$-CLC polynomials and Hurwitz-stable polynomials over $\K$. 
Let $f(t)=\sum_{k=0}^{d}a_{k}t^{k}$ be a real univariate polynomial of degree $d$. The even and odd parts of a real polynomial $f(t)$ are defined as:
  \begin{eqnarray*}
      f^{\even}(t)&:=a_0+a_2t^2+a_4t^4+\dots \\
      f^{\odd}(t)&:=a_1t+a_3t^3+a_5t^5+\dots 
  \end{eqnarray*}
  Define 
  \begin{eqnarray*}
      &f^{e}(t):=\sum_{m=0}^{\lfloor{ d/2 \rfloor}}(-1)^m a_{2m}t^m \\
      &f^{o}(t):=\sum_{m=0}^{\lfloor{ (d-1)/2 \rfloor}}(-1)^m a_{2m+1}t^{m} 
  \end{eqnarray*}
 \begin{theorem} \label{them:hermitebiehler}{Hermite-Biehler's criterion for stability:} \cite[Chap VII]{Levindistribution}, \cite{Kurtzsufficient}
      Let $f(t)=\sum_{k=0}^{d}a_{k}t^{k}$ be a univariate polynomial with positive $a_k$ for $k=0, \dots,d$. $f(t)$ is stable if and only if  $f^{e}(t)$ and $tf^{o}(t)$ have simple real interlacing zeros. Alternatively, $f(t)$ is stable if and only if all the zeros of $f^{\even}(t)$ and $f^{\odd}(t)$ are distinct, lie on the imaginary axis and alternate along it.
 \end{theorem} 
We use the following Hurwitz-Stability test for real polynomials based on the Interlacing Theorem and therefore, on the Boundary Crossing Theorem, cf. \cite[Chap 1]{Bhattacharyyarobust}.
 Note that $f(t)=f^{\even}(t)+f^{\odd}(t)$. Define the polynomial $g(t)$ of degree $d-1$ by :
 \begin{eqnarray*}
     \mbox{If} \ d=2m: g(t)&=[f^{\even}(t)-\frac{a_{2m}}{a_{2m-1}}tf^{\odd}(t)]+f^{\odd}(t)] \\
     \mbox{If} \ d=2m+1: g(t)&=[f^{\odd}(t)-\frac{a_{2m+1}}{a_{2m}}tf^{\even}(t)]+f^{\even}(t)]
 \end{eqnarray*}
 i.e., in general with $\mu =\frac{a_{d}}{a_{d-1}}$,
 $$g(t)=a_{d-1}t^{d-1}+(a_{d-2}-\mu a_{d-3})t^{d-2}+a_{d-3}t^{d-3}+(a_{d-4}-\mu a_{d-5})t^{d-4}+\dots$$
 \begin{theorem} \cite[Chap 1]{Bhattacharyyarobust} \label{them:stability}
   If $f(t)$ has all its coefficients positive, $f(t)$ is stable if and only if $g(t)$ is stable
 \end{theorem}
 It is well known and easy to verify that polynomials of degree $1$ and $2$ with positive leading coefficients are Hurwitz-stable if and only if all their coeﬃcients are positive. Thus the result is trivially true for $d=1$ and $d=2$. By \Cref{them:Routh-Hurwitz}, it's easy to verify that for $d=3$, the $3\times 3$ leading principal minor of $\Hur_{f(t)}$ is positive if and only if $a_1a_2-a_0a_3 >0$, cf. \cite[Chap-1]{Bhattacharyyarobust}. Since $a_1a_2 >3a_0a_3$ by \Cref{prop:logconcaveinequalities}, so $f(t)$ is Hurwitz-stable.

\begin{proposition} \label{prop:Hurwitzlowerdegree}
    Let $f(t)=\sum_{k=0}^{4}a_{k}t^{k}$ be a quartic polynomial with $a_k>0$ for $k=0,\dots,4$. 
    Suppose that
    \[
      a_{k}^2 \;>\; \frac{k+1}{k}\,a_{k-1}a_{k+1}
      \qquad\text{for }k=1,2,3.
    \]
    Then $f(t)$ is Hurwitz--stable.
\end{proposition}

\begin{proof}
By \Cref{them:stability}, $f(t)=\sum_{k=0}^{4}a_k t^k$ is Hurwitz--stable if and only if
\[
   g(t) \;:=\; f(t) - \mu t f^{\odd}(t)
   \;=\; a_0 + a_1 t + \Bigl(a_2 - \frac{a_4}{a_3} a_1\Bigr) t^2 + a_3 t^3
\]
is Hurwitz--stable, where $\mu = a_4/a_3$. Set
\[
  b_0 := a_0,\quad
  b_1 := a_1,\quad
  b_2 := a_2 - \frac{a_4}{a_3} a_1,\quad
  b_3 := a_3,
\]
so $g(t)=b_0 + b_1 t + b_2 t^2 + b_3 t^3$.

First, $b_2>0$ by \Cref{prop:logconcaveinequalities}(a), which gives
\[
   a_2 a_3 \;>\; a_1 a_4
   \quad\Longrightarrow\quad
   a_2 - \frac{a_4}{a_3} a_1 >0.
\]
Hence $b_0,b_1,b_2,b_3$ are all positive.

For a cubic with positive coefficients $b_0,b_1,b_2,b_3>0$, the Routh--Hurwitz criterion says that Hurwitz stability is equivalent to
\[
   b_1 b_2 \;>\; b_0 b_3
   \quad\Longleftrightarrow\quad
   a_1\Bigl(a_2 - \frac{a_4}{a_3}a_1\Bigr) \;>\; a_0 a_3.
\]
Equivalently,
\[
   a_1 a_2 a_3 - a_1^2 a_4 - a_0 a_3^2 \;>\; 0.
\]
Thus it remains to show
\begin{equation}\label{eq:quartic-Hurwitz-target}
   a_1 a_2 a_3 \;>\; a_1^2 a_4 + a_0 a_3^2.
\end{equation}

Define the dimensionless quantities
\[
   r_1 := \frac{a_1^2 a_4}{a_1 a_2 a_3} = \frac{a_1 a_4}{a_2 a_3},
   \qquad
   r_2 := \frac{a_0 a_3^2}{a_1 a_2 a_3} = \frac{a_0 a_3}{a_1 a_2}.
\]
Then \eqref{eq:quartic-Hurwitz-target} is equivalent to
\[
   1 - r_1 - r_2 \;>\; 0.
\]

From the assumed inequalities
\[
   a_1^2 > 2 a_0 a_2,\qquad
   a_2^2 > \frac{3}{2}a_1 a_3,\qquad
   a_3^2 > \frac{4}{3}a_2 a_4,
\]
we obtain the following bounds.

1. For $r_2$:
   \[
      \frac{a_0}{a_1} < \frac{a_1}{2a_2}
      \quad\Rightarrow\quad
      r_2 = \frac{a_0}{a_1}\cdot\frac{a_3}{a_2}
      < \frac{a_1}{2a_2}\cdot\frac{a_3}{a_2}
      = \frac{a_1 a_3}{2 a_2^2}.
   \]
   Using $a_2^2 > \frac{3}{2} a_1 a_3$ gives
   \[
      \frac{a_1 a_3}{a_2^2} < \frac{2}{3}
      \quad\Rightarrow\quad
      r_2 < \frac{1}{2}\cdot\frac{2}{3} = \frac{1}{3}.
   \]

2. For $r_1$:
   \[
      r_1 = \frac{a_1}{a_2}\cdot\frac{a_4}{a_3}.
   \]
   From $a_2^2 > \frac{3}{2}a_1 a_3$,
      $ \frac{a_1}{a_2} < \frac{2}{3}\cdot\frac{a_2}{a_3}$,
   and from $a_3^2 > \frac{4}{3}a_2 a_4$,
      $ \frac{a_4}{a_3} < \frac{3}{4}\cdot\frac{a_3}{a_2}$.
   Multiplying these estimates gives
   \[
      r_1 < \frac{2}{3}\cdot\frac{3}{4}\cdot\frac{a_2}{a_3}\cdot\frac{a_3}{a_2}
      = \frac{1}{2}.
   \]

Combining the two bounds, we get
\[
   r_1 + r_2 < \frac{1}{2} + \frac{1}{3} = \frac{5}{6} < 1,
\]
so $1 - r_1 - r_2 > 0$, which is exactly \eqref{eq:quartic-Hurwitz-target}. Hence $g$ is Hurwitz--stable, and therefore $f$ is Hurwitz--stable by \Cref{them:stability}.
\end{proof}

\begin{corollary}
    Let $f(t)=\sum_{k=0}^{d}a_{k}t^{k}$ be a univariate strictly CLC polynomial for $t \geq 0$ with $a_k>0$ for $k=0,\dots,d$. If $1 \leq d \leq 4$, then $f(t)$ is Hurwitz--stable.
\end{corollary}

\begin{lemma} \label{lemma:inequalityratio}
Assume $a_0,a_1,a_2,a_3,a_4,a_5 > 0$ and $a_2 - a_4 k > 0$. For any positive $k$,
\[
  \frac{a_1}{a_0} > \frac{a_3 - a_5 k}{a_2 - a_4 k}
  \quad\Longleftrightarrow\quad
  k < \frac{a_1 a_2 - a_0 a_3}{a_1 a_4 - a_0 a_5},
\]
provided $a_1 a_4 - a_0 a_5 > 0$.
\end{lemma}

\begin{proof}
By the standard identity for positive denominators,
\[
  \frac{a}{b} > \frac{c}{d} \;\Longleftrightarrow\; ad > bc.
\]
Applying this with $a = a_1$, $b = a_0$, $c = a_3 - a_5 k$, $d = a_2 - a_4 k$ and rearranging gives
\[
  a_1(a_2 - a_4 k) > a_0(a_3 - a_5 k)
  \;\Longleftrightarrow\;
  (a_1 a_4 - a_0 a_5)k < a_1 a_2 - a_0 a_3,
\]
which is equivalent to the claimed inequality when $a_1 a_4 - a_0 a_5 > 0$.
\end{proof}

\begin{theorem} \label{them:Hurwitz}
Let $f(t) = \sum_{k=0}^{5} a_k t^k$ be a univariate polynomial with $a_k > 0$ for $k = 0,\dots,5$. Suppose that
\[
  a_k^2 > \frac{k+1}{k}\,a_{k-1} a_{k+1}
  \quad \text{for } k = 1,\dots,4,
\]
and
\[
  (a_1 a_4 - a_0 a_5)^2 < (a_1 a_2 - a_0 a_3)(a_3 a_4 - a_2 a_5).
\]
Then $f(t)$ is Hurwitz-stable.
\end{theorem}

\begin{proof}
By \Cref{them:Routh-Hurwitz} and \Cref{remark:Routh-Hurwitz}, for $d=5$ it is sufficient to show that all leading principal minors $ \Delta_k > 0,\quad k=1,\dots,4,$ of the Hurwitz matrix
\[
  \Hur_{f(t)}
  =
  \begin{bmatrix}
    a_1 & a_3 & a_5 & 0   & 0 \\
    a_0 & a_2 & a_4 & 0   & 0 \\
    0   & a_1 & a_3 & a_5 & 0 \\
    0   & a_0 & a_2 & a_4 & 0 \\
    0   & 0   & a_1 & a_3 & a_5
  \end{bmatrix}.
\]

The Newton-type inequalities
  $a_k^2 > \frac{k+1}{k}\,a_{k-1} a_{k+1},\quad k=1,\dots,4,$ imply in particular that
\[
  a_1 > 0,\quad a_1 a_2 - a_0 a_3 > 0,\quad
  \det
  \begin{bmatrix}
    a_1 & a_3 & 0 \\
    a_0 & a_2 & a_4 \\
    0   & a_1 & a_3
  \end{bmatrix}
  > 0.
\]
From this it follows that $ \det
  \begin{bmatrix}
    a_1 & a_3 & a_5 \\
    a_0 & a_2 & a_4 \\
    0   & a_1 & a_3
  \end{bmatrix}
  > 0$, so $\Delta_k > 0$ for $k=1,2,3$.

It remains to show that the $4 \times 4$ leading principal minor is positive, i.e.,
\begin{equation} \label{eq:compodiv}
  \frac{a_1}{a_0}
  >
  \frac{
    \det \begin{bmatrix} a_3 & a_5 & 0 \\ a_1 & a_3 & a_5 \\ a_0 & a_2 & a_4 \end{bmatrix}
  }{
    \det \begin{bmatrix} a_2 & a_4 & 0 \\ a_1 & a_3 & a_5 \\ a_0 & a_2 & a_4 \end{bmatrix}
  }
  =
  \frac{a_3 - a_5 r}{a_2 - a_4 r},
\end{equation}
where $ r := \frac{a_1 a_4 - a_0 a_5}{a_3 a_4 - a_2 a_5}$.
By the Newton-type inequalities and \Cref{prop:logineq}(a),(c), we have
\[
  \frac{a_3 - a_5 r}{a_2 - a_4 r} > \frac{a_3}{a_2}.
\]
On the other hand, by \Cref{lemma:inequalityratio}, the inequality \eqref{eq:compodiv} holds if and only if
\[
  r < \frac{a_1 a_2 - a_0 a_3}{a_1 a_4 - a_0 a_5},
\]
that is,
\[
  \frac{a_1 a_4 - a_0 a_5}{a_3 a_4 - a_2 a_5}
  <
  \frac{a_1 a_2 - a_0 a_3}{a_1 a_4 - a_0 a_5}
  \;\Longleftrightarrow\;
  (a_1 a_4 - a_0 a_5)^2
  <
  (a_1 a_2 - a_0 a_3)(a_3 a_4 - a_2 a_5),
\]
which is exactly our additional hypothesis. Hence the $4 \times 4$ leading principal minor is positive and therefore $f(t)$ is Hurwitz-stable.
\end{proof}
In fact, a stronger result is presented in \cite[Theorem 3]{Hurwitzsuff} for entire functions. However, we limit our focus here to univariate polynomials, cf. \cite[Theorem 1]{Hurwitzsuff}.
\begin{theorem} \label{them:Hurwitzsuff}
If the coefficients of $f(t) = \sum_{k=0}^{d} a_k t^k$ are positive and satisfy
\[
  a_k a_{k+1} \;\geq\; \alpha\, a_{k-1} a_{k+2},\qquad 1 \leq k \leq d-1,
\]
where $\alpha$ is the unique positive root of the cubic
\[
  t^3 - t^2 - 2t - 1 = 0 \quad (\alpha \approx 2.1479),
\]
then $f(t)$ is Hurwitz-stable. In particular, the conclusion holds if
\[
  a_k^2 \;\geq\; \sqrt{\alpha}\, a_{k-1} a_{k+1},\qquad 1 \leq k \leq d-1.
\]
\end{theorem}

\noindent

Now we extend the above univariate results to multivariate homogeneous polynomials.

\begin{theorem} \label{them:Hurwitzclc}
  Let $f(x) \in \R^{d}_{n}[x]$ be a strictly $\K$-Lorentzian polynomial of degree
  $d$ with $1 \le d \le 4$ over a proper convex cone $\K$. Then $f(x)$ is
  Hurwitz-stable over $\K$.
\end{theorem}

\begin{proof}
Since $f$ is strictly $\K$-Lorentzian, \Cref{def:clcline} together with
\Cref{them:plsclc} (the equivalence of strictly $\K$-CLC and strictly
$\K$-Lorentzian polynomials in the homogeneous setting) implies that for any
$x,v \in \inter \K$ the univariate restriction
\[
  f(x+tv)
  = \sum_{k=0}^{d} \frac{D_v^k f(x)}{k!}\, t^k
\]
is a strictly Lorentzian polynomial in $t$ for all $t$ such that $x+tv \in \K$.
Moreover, by \Cref{cor:noncoeffs}, all coefficients
\(
  a_k := \frac{D_v^k f(x)}{k!}
\)
are positive for $x,v \in \inter \K$.

Since $f$ is strictly $\K$-Lorentzian, \Cref{thm:logconcave} yields the
Newton-type inequalities
\[
  \bigl(D_v^k f(x)\bigr)^2
  > \frac{k+1}{k}\, D_v^{k-1} f(x)\, D_v^{k+1} f(x),
  \qquad 1 \le k \le d-1,
\]
for all $x,v \in \inter \K$. In terms of the coefficients
$a_k = \frac{D_v^k f(x)}{k!}$, this is exactly
\[
  a_k^2 > \frac{k+1}{k}\, a_{k-1} a_{k+1},\qquad 1 \le k \le d-1.
\]
Therefore, by \Cref{prop:Hurwitzlowerdegree}, the univariate polynomial
$t \mapsto f(x+tv)$ is Hurwitz-stable for every $x,v \in \inter \K$.
By \Cref{def:hurstable}, this means that $f(x)$ is Hurwitz-stable over $\K$.
\end{proof}

\begin{theorem} \label{them:hurwitzclc5}
Let $f(x) \in \R^{d}_{n}[x]$ be strictly $\K$-Lorentzian of degree $d \ge 5$
over a proper convex cone $\K$. Suppose that
\[
  \bigl(D_v^{k} f(x)\bigr)^2
  > \frac{\sqrt{\alpha}\, k}{k+1}\,
    D_v^{k-1} f(x)\, D_v^{k+1} f(x)
  \quad \text{for all } 1 \le k \le d-1,\; x,v \in \inter \K,
\]
where $\alpha$ is the unique positive root of the cubic
$t^3 - t^2 - 2t - 1 = 0$ (so $\alpha \approx 2.1479$). Then $f(x)$ is
Hurwitz-stable over $\K$.
\end{theorem}

\begin{proof}
Fix $x,v \in \inter \K$ and consider the univariate restriction
\[
  g(t) := f(x+tv)
  = \sum_{k=0}^{d} a_k t^k,
  \qquad \text{with } a_k := \frac{D_v^k f(x)}{k!}.
\]
By \Cref{def:clcline} and \Cref{them:plsclc}, $g$ is strictly Lorentzian, so in
particular $a_k > 0$ for all $k$.

The assumed inequalities on $D_v^k f(x)$ translate into the coefficient
inequalities required in \Cref{them:Hurwitzsuff}. Indeed,
\[
  a_k^2
  = \frac{\bigl(D_v^k f(x)\bigr)^2}{(k!)^2}
  > \sqrt{\alpha}\,\frac{k}{k+1}
     \frac{D_v^{k-1} f(x)\, D_v^{k+1} f(x)}{(k!)^2}
  =
  \sqrt{\alpha}\,
  \frac{D_v^{k-1} f(x)}{(k-1)!}\,
  \frac{D_v^{k+1} f(x)}{(k+1)!}
  = \sqrt{\alpha}\, a_{k-1} a_{k+1}.
\]
Thus the sufficient condition in \Cref{them:Hurwitzsuff} is satisfied, and so
$g(t)$ is Hurwitz-stable. Since this holds for every $x,v \in \inter \K$, it
follows from \Cref{def:hurstable} that $f(x)$ is Hurwitz-stable over $\K$.
\end{proof}

\begin{remark}
Not every strictly $\K$-CLC polynomial is Hurwitz-stable over $\K$.
For instance, for $\K = \R_{\ge 0}$, the univariate polynomial
\[
  f(t) = 5 + 14 t + 12.5 t^2 + 7.2 t^3 + 3 t^4 + t^5
\]
is Lorentzian for $t \ge 0$ but not Hurwitz-stable. Conversely, the polynomial
\[
  f(t) = 5 + 25 t + 50 t^2 + 30 t^3 + 10 t^4 + 3 t^5
\]
is Hurwitz-stable but not Lorentzian. These examples show that the sufficient
conditions in \Cref{them:Hurwitzclc,them:hurwitzclc5} are not necessary.
\end{remark}

\section{Applications: Stability Analysis for LTI Systems} \label{sec:dynamics}

Consider the linear time-invariant (LTI) system
\[
  \dot{x} = A x, \qquad A \in \R^{n \times n}.
\]
A matrix $A$ is called \emph{stable} (or \emph{Hurwitz}) if every eigenvalue of $A$
has negative real part, i.e., all eigenvalues lie in the open left half-plane.
If $A$ is real symmetric, then $A$ is stable if and only if $A$ is negative
definite. Lyapunov's theorem asserts that $A$ is stable if and only if there
exists a symmetric positive definite matrix $P$ such that the Lyapunov operator
\[
  L_A(P) := A^\top P + P A
\]
is negative definite.

\begin{theorem}[Routh--Hurwitz criterion] \label{them:Routh-Hurwitz-stable}
Let $A \in \R^{n \times n}$, and let $\chi_A(t)$ be its characteristic
polynomial. Then the following are equivalent \cite[Chap.~9]{Poznyak},
\cite[Chap.~7]{stabilityHurwitz}:
\begin{enumerate}[(a)]
\item The matrix $A$ is stable, i.e., all eigenvalues of $A$ lie in the open
      left half-plane.
\item The characteristic polynomial $\chi_A(t)$ is Hurwitz-stable.
\item All the leading principal minors of the Hurwitz matrix $\Hur_{\chi_A(t)}$
      are positive.
\end{enumerate}
\end{theorem}
Hence we have the following.
\begin{theorem} \label{them:lconstable}
Let $f(t)$ be a strictly CLC polynomial of degree $d \leq 4$ for $t \geq 0$.  Then the corresponding dynamical system $\frac{dx}{dt}=Ax$ is stable if the characteristic polynomial $\chi(A)=f(t)$. 
\end{theorem}
\begin{proof}
    Since $f(t)$ is a strictly CLC polynomial of degree $d \leq 4$ for $t \geq 0$, $f(t)$ is Hurwitz-stable by \Cref{them:Hurwitzclc}. Then it follows from \Cref{them:Routh-Hurwitz-stable}.
\end{proof}
\begin{corollary}
Let $f(t)=\sum_{k=0}^{5} a_{k} t^{k}$ be a strictly CLC polynomial of degree $5$
for $t \geq 0$. If
\[
  (a_1 a_4 - a_0 a_5)^2
  < (a_1 a_2 - a_0 a_3)\,(a_3 a_4 - a_2 a_5),
\]
and $\chi_A(t) = f(t)$ is the characteristic polynomial of a matrix $A$, then
the LTI system $\dot{x} = A x$ is stable.
\end{corollary}

\begin{proof}
By Theorem~\ref{them:Hurwitz}, the inequality implies that $f(t)$ is
Hurwitz-stable. Hence $\chi_A(t)$ is Hurwitz-stable, and the claim follows from
Theorem~\ref{them:Routh-Hurwitz-stable}.
\end{proof}
\begin{theorem} \label{them:lconstable5}
Let $f(t)=\sum_{k=0}^{d}a_{k}t^{k}, d \geq 5$ be a univariate polynomial with positive coefficients and $a_k^2 > \sqrt{\alpha} a_{k-1}a_{k+1}$ where $\alpha$ is the unique positive root of the polynomial $t^3- t^2- 2t-1 (\alpha \approx 2.1479)$. 
Then the corresponding dynamical system $\frac{dx}{dt}=Ax$ is stable if the characteristic polynomial $\chi(A)=f(t)$,.
\end{theorem}
\begin{proof}
It follows directly from Theorem~\ref{them:Hurwitzsuff} that $f(t)$ is Hurwitz-stable since the
assumption $a_k^2 > \sqrt{\alpha}\, a_{k-1} a_{k+1}$ for all $1 \le k \le d-1$
implies the sufficient condition in that theorem.  Then it follows from
Theorem~\ref{them:Routh-Hurwitz-stable}.
\end{proof}
\noindent Furthermore, we extend the result for multivariate polynomials.
\begin{theorem}
Let $f(x)$ be a strictly $\K$-Lorentzian polynomial of degree $d \leq 4$ over a
proper convex cone $\K$. Suppose there exist $x_0, v \in \inter \K$ such that
the characteristic polynomial of a matrix $A$ is given by
\[
  \chi_A(t) = f(x_0 + t v).
\]
Then the LTI system $\dot{x} = A x$ is stable.
\end{theorem}

\begin{proof}
By Theorem~\ref{them:Hurwitzclc}, $f$ is Hurwitz-stable over $\K$, so for any
$x,v \in \inter \K$ the univariate restriction $t \mapsto f(x+tv)$ is
Hurwitz-stable. In particular, $\chi_A(t) = f(x_0+tv)$ is Hurwitz-stable.
The conclusion then follows from
Theorem~\ref{them:Routh-Hurwitz-stable}.
\end{proof}
\begin{theorem}
Let $f(x)$ be a strictly $\K$-Lorentzian polynomial of degree $d \ge 5$ over a
proper convex cone $\K$. Suppose there exist $x_0, v \in \inter \K$ such that
\[
  \chi_A(t) = f(x_0 + t v) = \sum_{k=0}^{d} a_k t^k
\]
is the characteristic polynomial of a matrix $A$, and the coefficients satisfy
\[
  a_k^2 > \sqrt{\alpha}\, a_{k-1} a_{k+1}, \qquad 1 \le k \le d-1,
\]
where $\alpha$ is the unique positive root of $t^3 - t^2 - 2t - 1 = 0$
($\alpha \approx 2.1479$). Then the LTI system $\dot{x} = A x$ is stable.
\end{theorem}

\begin{proof}
By the coefficient inequalities and Theorem~\ref{them:Hurwitzsuff},
$\chi_A(t)$ is Hurwitz-stable. Stability of $\dot{x} = A x$ then follows from
Theorem~\ref{them:Routh-Hurwitz-stable}.
\end{proof}
\begin{remark}
Whenever the characteristic polynomial $\chi_A(t)$ of a stable matrix $A$ can
be realized as a univariate restriction $\chi_A(t) = f(x_0 + t v)$ of a
Hurwitz-stable polynomial $f$ over a proper convex cone $\K$ with
$x_0, v \in \inter \K$, we may regard $f$ as a $\K$-Lorentzian lift of the
dynamical system $\dot{x} = A x$. This viewpoint connects LTI stability with
the cone-geometric properties developed in this paper.
\end{remark}
In particular, for a linear time-invariant system $\dot{x} = A x$, Hurwitz stability of the characteristic polynomial $\chi_A(t)$ is equivalent to asymptotic stability of the equilibrium $x = 0$. Thus, our Newton-type coefficient inequalities for strictly CLC and strictly $\K$-Lorentzian polynomials yield simple sufficient tests ensuring stability of LTI systems directly from the coefficients of $\chi_A(t)$.   Conversely, given a strictly $\K$-CLC polynomial with suitable Hurwitz-stability properties, we can construct LTI systems whose spectral behaviour is encoded by $f$.

\noindent \textbf{Acknowledgment}:
During the preparation of this manuscript, the author used ChatGPT to improve the clarity and conciseness of the writing. All content was subsequently reviewed, revised, and verified by the author, who takes full responsibility for the final version of the article.

\bibliographystyle{alpha}
\bibliography{main}

\begin{thebibliography}{ALOGV24}

\bibitem[AGV21]{Cynthialog}
Nima Anari, Shayan~Oveis Gharan, and Cynthia Vinzant.
\newblock Log-concave polynomials, {I}: entropy and a deterministic approximation algorithm for counting bases of matroids.
\newblock {\em Duke Math. J.}, 170(16):3459--3504, 2021.

\bibitem[ALOGV24]{Cynthialog3}
Nima Anari, Kuikui Liu, Shayan Oveis~Gharan, and Cynthia Vinzant.
\newblock Log-concave polynomials {III}: {M}ason's ultra-log-concavity conjecture for independent sets of matroids.
\newblock {\em Proc. Amer. Math. Soc.}, 152(5):1969--1981, 2024.

\bibitem[BD24]{GPlorentzian}
Grigoriy Blekherman and Papri Dey.
\newblock $\mathcal{K}$-lorentzian polynomials.
\newblock {\em arXiv preprint arXiv:2405.12973}, 2024.

\bibitem[BF76]{Barkerselfdual}
George~Philip Barker and James Foran.
\newblock Self-dual cones in {E}uclidean spaces.
\newblock {\em Linear Algebra Appl.}, 13(1-2):147--155, 1976.
\newblock Collection of articles dedicated to Olga Taussky Todd.

\bibitem[BGLS01]{Gulerhyperbolicconvex}
Heinz~H Bauschke, Osman G{\"u}ler, Adrian~S Lewis, and Hristo~S Sendov.
\newblock Hyperbolic polynomials and convex analysis.
\newblock {\em Canadian Journal of Mathematics}, 53(3):470--488, 2001.

\bibitem[BH20]{Brandenlorentzian}
Petter Br{\"a}nd{\'e}n and June Huh.
\newblock Lorentzian polynomials.
\newblock {\em Annals of Mathematics}, 192(3):821--891, 2020.

\bibitem[BK64]{Bellmanselected}
Richard~Ernest Bellman and Robert~E Kalaba.
\newblock Selected papers on mathematical trends in control theory.
\newblock {\em (No Title)}, 1964.

\bibitem[BK95]{Bhattacharyyarobust}
Shankar~P Bhattacharyya and Lee~H Keel.
\newblock Robust control: the parametric approach.
\newblock In {\em Advances in control education 1994}, pages 49--52. Elsevier, 1995.

\bibitem[BL23a]{Leake}
Petter Br{\"a}nd{\'e}n and Jonathan Leake.
\newblock Lorentzian polynomials on cones.
\newblock {\em arXiv preprint arXiv:2304.13203}, 2023.

\bibitem[BL23b]{Hereditary}
Petter Br{\"a}nd{\'e}n and Jonathan Leake.
\newblock Lorentzian polynomials on cones.
\newblock {\em arXiv preprint arXiv:2304.13203}, 2023.

\bibitem[BP94]{Bermannonnegative}
Abraham Berman and Robert~J. Plemmons.
\newblock {\em Nonnegative matrices in the mathematical sciences}, volume~9 of {\em Classics in Applied Mathematics}.
\newblock Society for Industrial and Applied Mathematics (SIAM), Philadelphia, PA, 1994.
\newblock Revised reprint of the 1979 original.

\bibitem[BS75]{Barkeralgebraic}
GP~Barker and Hans Schneider.
\newblock Algebraic perron-frobenius theory.
\newblock {\em Linear Algebra and its Applications}, 11(3):219--233, 1975.

\bibitem[BV04]{Boydconvex}
Stephen~P Boyd and Lieven Vandenberghe.
\newblock {\em Convex optimization}.
\newblock Cambridge university press, 2004.

\bibitem[COSW04]{Oxley}
Young-Bin Choe, James~G Oxley, Alan~D Sokal, and David~G Wagner.
\newblock Homogeneous multivariate polynomials with the half-plane property.
\newblock {\em Advances in Applied Mathematics}, 32(1-2):88--187, 2004.

\bibitem[Gan98]{gantmacher}
F.~R. Gantmacher.
\newblock {\em The theory of matrices. {V}ol. 1}.
\newblock AMS Chelsea Publishing, Providence, RI, 1998.
\newblock Translated from the Russian by K. A. Hirsch, Reprint of the 1959 translation.

\bibitem[GB04]{stabilitylevi}
Daniel Goeleven and Bernard Brogliato.
\newblock Stability and instability matrices for linear evolution variational inequalities.
\newblock {\em IEEE Transactions on Automatic Control}, 49(4):521--534, 2004.

\bibitem[Gre81]{Gregor}
Ji\v{r}\'{\i} Gregor.
\newblock On quadratic {H}urwitz forms. {I}.
\newblock {\em Apl. Mat.}, 26(2):142--153, 1981.
\newblock With a loose Russian summary.

\bibitem[Gur09]{Gurvitsmultivariate}
Leonid Gurvits.
\newblock On multivariate newton-like inequalities.
\newblock In {\em Advances in combinatorial mathematics}, pages 61--78. Springer, 2009.

\bibitem[HH69]{Remarkcopositive}
Emilie Haynsworth and A.~J. Hoffman.
\newblock Two remarks on copositive matrices.
\newblock {\em Linear Algebra Appl.}, 2:387--392, 1969.

\bibitem[JM12]{stabilityHurwitz}
Rolf Jeltsch and Mohamed Mansour.
\newblock {\em Stability Theory: Hurwitz Centenary Conference Centro Stefano Franscini, Ascona, 1995}, volume 121.
\newblock Birkh{\"a}user, 2012.

\bibitem[KR48]{Kreinlinear}
Mark~Grigor’evich Krein and Mark~A Rutman.
\newblock Linear operators leaving invariant a cone in a banach space.
\newblock {\em Uspekhi mat. nauk}, 3(1):3--95, 1948.

\bibitem[Kur92]{Kurtzsufficient}
David~C. Kurtz.
\newblock A sufficient condition for all the roots of a polynomial to be real.
\newblock {\em Amer. Math. Monthly}, 99(3):259--263, 1992.

\bibitem[KV08]{Hurwitzsuff}
Olga~M. Katkova and Anna~M. Vishnyakova.
\newblock A sufficient condition for a polynomial to be stable.
\newblock {\em J. Math. Anal. Appl.}, 347(1):81--89, 2008.

\bibitem[Lev80]{Levindistribution}
B.~Ja. Levin.
\newblock {\em Distribution of zeros of entire functions}, volume~5 of {\em Translations of Mathematical Monographs}.
\newblock American Mathematical Society, Providence, RI, revised edition, 1980.
\newblock Translated from the Russian by R. P. Boas, J. M. Danskin, F. M. Goodspeed, J. Korevaar, A. L. Shields and H. P. Thielman.

\bibitem[Poz09]{Poznyak}
Alexander~S Poznyak.
\newblock {\em Advanced mathematical tools for automatic control engineers: Stochastic techniques}.
\newblock Elsevier Ltd, 2009.

\bibitem[SV70]{Cross}
Hans Schneider and Mathukumalli Vidyasagar.
\newblock Cross-positive matrices.
\newblock {\em SIAM Journal on Numerical Analysis}, 7(4):508--519, 1970.

\bibitem[Tar18]{SymmetricPFstructure}
Pablo Tarazaga.
\newblock On the structure of the set of symmetric matrices with the perron--frobenius property.
\newblock {\em Linear Algebra and its Applications}, 549:219--232, 2018.

\bibitem[Van68]{Vandergraftspectral}
James~S. Vandergraft.
\newblock Spectral properties of matrices which have invariant cones.
\newblock {\em SIAM J. Appl. Math.}, 16:1208--1222, 1968.

\end{thebibliography}

\end{document}